\documentclass{amsart}

\title{AP Theory II: Intrinsic 4D Quantum YM Theory with Mass Gap}
\author{H. E. Winkelnkemper}

\address{\begin{flushleft} \quad Department of Mathematics\\
\quad University of Maryland\\
\quad College Park, Maryland 20742
\end{flushleft}}

\email{hew@math.umd.edu}

\begin{document}
\maketitle

\begin{abstract} We describe a sub-theory of Artin Presentation Theory (AP Theory), which has many genuine,discrete,group-theoretic,non-infinitesimal, {\it qualitative} analogues (including with the mass gap) of the main desiderata of the hypothetical {\it quantitative} infinitesimal '4D Quantum YM Theory' for the so-called Clay Millenium 'YM Existence and Mass Gap' problem.
Our entirely mathematically rigorous theory is {\it not a model},no new axioms or measures are introduced,does not rely on SUSY,is free of smooth 4D singularities,moduli spaces,path integrals, graph/lattice combinatorics and probabilistic,category,twistor or topos arguments and is intimately related to the theory of pure framed braids.Despite being based on a rigorous, radical,universal Holographic Principle,the theory still contains an analogue of Donaldson/Seiberg-Witten Theory, an {\it infinitely generated}, at each stage, graded group of topology-changing transitions and interactions and more. The theory is a {\it sporadic} algebraic theory graded by the positive integers, characteristic of dimension 4, non-perturbative, background independent, parameter-free and so basic and conceptually simple,that it cannot be reduced any further in a rigorous
fashion: just discrete, finitely presented group theory, augmented by a non-infinitesimal, non-local, 'exterior' construction of 4D {\it smooth} structures on compact, simply-connected 4-manifolds, with a connected boundary. This radical conceptual simplicity immediately avoids ultraviolet divergence and renormalization problems. It thus becomes natural to
conjecture that any other such theory,which purports to be a rigorous solution to the YM Millenium problem,has to interact in a non-trivial way with our topological gauge theory,
even if the latter is, a priori, 'exterior' and {\it intrinsic} (in the sense that Cobordism Theory is exterior and  intrinsic) and in this way sets the stage for a rigorous study of the fundamental problem of non-perturbative dynamics of classical 4D gauge theories.Our main contention is: the radical,universal AP-holography, with its strong topology changing interactions, which reach all the way to the 'vacuum' of discrete group theory,may destroy any infinitesimal,PDE based
Field Theory,required for solving the actual Clay YM problem in its present {\it quantitative} form as a problem of so-called 'constructive' 3+1 QFT. More generally,due to the fact that AP Theory is not a model,e.g., does not introduce any new axioms,any rigorous axiomatic $3+1$ QFT has to align itself with it in a mathematically rigorous fashion.
\end{abstract}

\section{Introduction}
In this somewhat long introduction,we start from general analogies between the theory of Artin Presentations (AP Theory) with modern physics, before giving more specific ones related to the actual Clay YM problem,\cite {JW}, \cite{D}, \cite{V}, \cite{F},in section 3. 

In the following $X^4$ will denote a connected,compact,smooth,simply-connected 4-manifold with a connected boundary,$\partial X^4$; if $\partial X^4$ = $S^3$,we will also denote the so determined closed manifold by $X^4$.Although all the $X^4$ we discuss are so-called '2-handle bodies' ,(see, e.g., \cite{GS})
they form a very large class of smooth 4-manifolds: every connected, closed,orientable 3-manifold 
appears as a $\partial X^4$ and at least every complex elliptic surface (e.g.the Kummer surface) is 
part of this class, \cite{CW}.(It is unknown whether every {\it closed}, smooth,simply-connected 4-manifold can be so obtained, although it seems very likely that this is so).\cite{GS},p.344.

In AP Theory (see ahead) such a $X^4$ is already determined by a certain type of (group) presentation, an Artin Presentation, $r$, on $n$ generators with $n$ relations, of the fundamental group of its boundary,which we denote by $\pi(r)$ and 
we write, $X^4$ = $W^4(r)$.

Even at this early stage,this relation of pure discrete group theory with 4D smooth structures has physical relevance:

It is the most radical,mathematically rigorous,universal form of $3+1$ holography:
Discrete Group theory on the boundary of $X^4$ determines the whole smooth 4-manifold $X^4$,up to diffeomorphism.

For the meaning and importance in physics of holography ('t Hooft, Susskind, Maldacena,\dots,) see \cite{Bo}.

The discrete,purely group-theoretic Artin presentation $r$ is a {\it hologram} of the {\it  smooth} 4D spacetime universe $W^4(r)$.Compare to \cite{M2},p.63.

This is the more far reaching AP-analogue of Witten's most general version of  holography in 
modern physics:

{\it "Gravity on the bulk is built from Gauge Theory on the boundary"},see p.13
of \cite{Wi1}.

{\it This analogy of the rigorous, radical, universal AP-holography with,e.g.the more heuristic, restrictive, AdS/CFT holography,(\cite{M1},\cite{M2}, \cite{Re},\cite{Wi2}, \cite{SW}, \cite{S1},\dots), is the prototype of the analogies exhibited and studied in this paper.}

We can say: the compact,smooth $3+1$ spacetime universe $W^4(r)$ {\it emerges} from the discrete group-theoretic 'vacuum', via the purely group theoretic Artin presentation $r$; in other words,4D gravity,in its most abstract Einsteinian sense,emerges purely group-theoretically,{\it not} metrically,and it does so in a non-infinitesimal,non-local manner.Compare to \cite{AJL},\cite {Lo},\cite{Se},\cite{Wi1},p.5,\cite{Sm}, \cite{M4}.

The 4D smooth {\it 'continuum picture arises from this fundamental discreteness'}, \cite {As}, p.174.

This {\it 'immense gauge symmetry'} (\cite{G},p.A98) already makes AP Theory a candidate for
 a topological  {\it 'quantum gauge theory in four dimensional space-time'}, \cite{JW},p.3, although it is, apriori, 'exterior',not confined to the Differential Geometry of a single fixed  manifold,in the same sense
that Cobordism Theory is not so confined. In AP Theory,gauge theory has exterior,autonomous, {\it intrinsic}, absolute meaning,not related to any particular Lie group ,just as Cobordism Theory has exterior,intrinsic,absolute homological meaning.We note here that,nevertheless, cobordism theory has led to the solutions of important 'internal' problems on a given fixed manifold (Hirzebruch formula,Atiyah-Singer theorems,
Steenrod's conjecture,etc. \cite {A1}, \cite{Su}).

As has been pointed out,even the more classical AdS/CFT holographic principle {\it 'calls into question not only the status of field theory but the
very notion of locality'}, \cite{Bo},p.2, \cite{Wi3},p.1579, \cite{S1},p.42, \cite{S2}, p.10.  Thus infinitesimal, analytic, smooth continuum using differential geometric methods are not, apriori, used in AP Theory, {\it although} the mathematical results of such classic analytic gauge theoretic methods,e.g. Donaldson/Seiberg-Witten Theory,can still lead to purely group-theoretic analogs in AP theory, {\it despite} the absence of moduli spaces in AP Theory.  See \cite {W1}, p.240, \cite{R}, p.621 and section 2 ahead.

In other words,the non-local radical AP-holography forces the substitution of Field Theories by discrete group theory. We do not consider this a bad thing, as, e.g. the ultraviolet catastrophe in QFT:to AP  discrete group theory it does not matter 'whether a field at a point is not well-defined' and has to be remedied 'by smearing over test functions to tame the UV divergences'....On the contrary, AP theory creates our intrinsic gauge theory,
without obtaining the gauge theory analytically  from vector bundles/fields,connections,on a particular given manifold. In the 4D AP Theory no field redefinition or renormalization is needed.
(compare to \cite{Wi2},p.2). Furthermore, our non-local construction of 4D smooth structures is still so subtle and metamathematically "local" so to speak, that, as mentioned above, an analogue of Donaldson/Seiberg-Witten Theory is still obtained.  This AP-theoretic non-locality is evidently the conceptually  
simplest way to obtain universal, {\it intrinsic quantum} gauge theory, not based on a particular Lie group.  Thus AP Theory takes the place of {\it "the large N limit of YM theory"}, \cite{M3}, \cite{M4}.

{\it In AP Theory the UV problem and other quantum field-theoretic problems are bypassed rigorously with discrete group theory}. (Compare to \cite{Fe}, \cite{P}).

In a more abstract vein, our metamathematical maximalization of group theory, which does not use infinite-dimensional spaces (i.e. does not rely on SUSY), but which nevertheless still contains a very general, at each stage, {\it infinitely generated} graded group of topology-changing transitions (the
AP-analogue of Morse Theory) relates the symmetries of our 'particles',i.e. the discrete Artin Presentations $r$, with the exterior symmetries of the smooth 4D spacetime universes $W^4(r)$.
This is an analogue of the very desired property for unifying the symmetries of the actual Standard Model with those of classical General Relativity.  For a promising relation between braids and the Standard Model, see \cite{BT}.

One can say, in AP Theory, SUSY is already 'broken' metamathematically by its canonical,natural grading by the positive 
integers; SUSY reliance on an infinite number of dimensions is substituted by the infinite number of generators of our group of topology-changing transitions, whose existence is requiered for any sensible  Quantum Gravity theory, see \cite{Wi2}, p.4.  There is no incompatibility problem between SUSY and the gauge transitions and other symmetries of AP Theory.(\cite{Wi5}, p.5, \cite{Wi6}, p.362, \cite{Wi18}).  The 'vacuum' of AP Theory, i.e. the purely discrete part, gets enough 'rigidity' from its own graded group theory. A priori, SUSY is merely the stabilization theory of AP Theory,when $n$ approaches infinity. 

These topology-changing transitions, which we call Torelli transitions, (and which,                                                                                                                                                                                                                                               for each $n$, form a group, isomorphic to the commutator subgroup of the pure braid group on $n$ strands, compare to \cite{N}, p.43) are more universal, versatile and sharper than those caused by any classical Morse Theory or $3+1$ TQFT, \cite {A2}, in the sense that, e.g. they can also just leave the underlying {\it topological} structure of the manifold $X^4$ invariant and just change the smooth structure or in fact also leave it fixed; (see section 2 ahead). It is discrete pure group theory, that in a general, systematic, exterior manner (thus avoiding symmetry destroying 'skein' and/or ad hoc 'by hand' internal surgery methods) generates new smooth 4D structures,i.e.structures that are at the foundations of modern 4D gravitational physics,\cite{C1},\cite{C2}.

In order to further explain and augment the above, it is instructive to analize more explicitly how the {\it smooth} 4D manifold $W^4(r)$ is obtained holographically in a non-infinitesimal, non-local manner from the {\it purely discrete} Artin presentation $r$ of the
fundamental group of its boundary.

Let $\Omega_n$ denote the compact 2-disk with $n$ holes in the plane and $\partial \Omega_n$ its boundary.  To obtain $W^4(r)$ from $r$ proceed as follows:an Artin presentation $r$ on $n$ generators defines a unique, {\it framed}, pure braid on $n$ strands (and conversely).This framing,(i.e. an assignement of an integer to each strand,sometimes also called a 'coloring' of the braid), is not obtained 'by hand', but is obtained canonically by representing the pure,framed, braid uniquely by an Artin presentation $r$, which then defines a smooth diffeomorphism $h(r):\Omega_n\to \Omega_n$,
which restricts to the identity on $\partial \Omega_n$. With $h(r)$,the smooth 4-manifold $W^4(r)$ is constructed by means of a relative open book construction (see \cite{W1}, p.250,\cite{CW}, \cite{C1},\cite{C2},and references therein, for the rigorous technical details).  This is a construction which is structurally similar to the fundamental Lefschetz Hyperplane Theorem for non-singular complex algebraic varieties and goes beyond the mere surgery prescriptions of the Kirby Calculus, see section 2 ahead.

{\it The diffeomorphism $h(r):\Omega_n\to \Omega_n$ is determined by $r$ only up to isotopies keeping it fixed as the identity on the boundary $\partial \Omega_n$.}

Thus (the still well-defined) smooth 4D diffeomorphism class of the 4-manifold $W^4(r)$ is obtained in a non-infinitesimal,non-local manner from the 2D diffeomorphism $h(r)$, which up to isotopy is determined by the discrete Artin presentation $r$.

This is a {\it sporadic} 4D smooth analogue, in a metamathematical sense, to the celebrated Hilbert Vth Problem, which cleared up Lie Group Theory's 'infinitesimal mess' (\cite{H}) by showing that the existence of the {\it smooth} structure of a Lie group need not be postulated a priori, since its existence already follows from non-infinitesimal, non-analytic concepts.

{\it This is the fundamental smooth 4D topological construction of AP Theory.}

It is indeed universal $3+1$ {\it'gravitational holography'},very different from the very restricted field theoretic version of
Maldacena ,Rehren,et al. and is definitely not {\it 'a metaphoric illusion'}. (\cite{S2},p.9).

As far as gauge theory is concerned,it is in the reverse order of the dimensional reduction of \cite{Be},which is a cornerstone of so-called Geometric Langlands Theory,\cite{GW}, \cite{W4}.

It is natural to conjecture that,just as in the mathematically analogous case of Lie Groups,our
sporadic holographic 4D gravitational Hilbert Vth-like construction,where 4D smooth structures emerge from pure discrete group theory,will help clear up the current conceptual mess in Quantum Gravity,(see \cite{Sm}, \cite{N}) beyond the results of this paper.

In order to place our main section  3 ahead in its proper metamathematical setting,we point out some more similarities of our fundamental construction to more heuristic ideas and concepts of modern physics in the  literature,(\cite{Wi7}, \cite{Wi8}, \cite{A2}, \cite{G1}, \cite{G2}, \cite{GW},\cite{Wi9},\dots )

i) In the above construction of the 4D smooth manifold $W^4(r)$ from the 2D diffeomorphism  $h(r):\Omega_n\to \Omega_n$,the boundary components of $\partial \Omega_n$ define $n+1$ knots in the 3D  boundary of $W^4(r)$ on which the,in general topology-changing and knot-changing, Torelli transitions/interactions  act.\cite{W1},p.226.

This is a  more direct, non-categorical analogue of {\it 'surface operators...realizing knot invariants'} of the 4D topological gauge theories of Gukov,Witten, et al., \cite{G1}, \cite{GW}, \cite{Wi9}, p.5, which lead to Geometric Langlands Theory.  The following  quote of p.9 of \cite{G1} also seems relevant here:{\it "..every topological gauge theory,which admits surface operators is,in a sense,a factory that produces examples of braid group actions on branes.."}.

ii) The Artin presentation $r$ on $n$ generators, which is obtained by infinite intersection via the {\it rigorous} Cayley-von Dyck process (see section 2), instead of the infinite 'sum' (as in the non-rigorous Feynman process for achieving independence of the metric),defines a pure framed   (i.e. colored) braid on $n$ strands,i.e.intuitively a back-ground independent, {\it macroscopic}  'string',{\it which is immediately related (without relying on SUSY), to 4D 'gravity' via the 4D smooth structure of $W^4(r)$, by the fundamental Hilbert Vth-like construction of AP Theory}, compare to \cite{Wi11}, p.25, \cite{Wi16}. These are the strings in AP Theory, when it is considered as {\it 'the infinite limit of $SU(N)$ YM Theory'},\cite{M3},p.10, \cite{M4}, \cite{M5}.

(A priori,they do not seem to be related to so-called 'string topology',\cite{Su1}.)

Thus AP Theory can be considered to be a graded,  background independent, non-perturbative, parameter-free, macroscopic  'string theory' where holography and topology changing transistions and interactions are as strong as possible, (compare to p.411 of \cite{Gr} and p.8 of \cite{We}). The intrinsic 'world-sheet' of the string $r$ consists  simply of the iterates of the canonically associated planar diffeomorphisms $h(r):\Omega_n\to \Omega_n$ .The incredibly rich Iteration and  Covering theories (Nielsen, Thurston,...) of these  $h(r):\Omega_n\to \Omega_n$  makes this 'string theory' a very strong one indeed; in fact,in should also be relevant to 
Loop Quantum Gravity,\cite{Sm}, \cite{N}, \cite{P},\cite{BT}, thus uniting these two,a priori different approaches of modern physics to Quantum Gravity. The true LQG,after it accomodates universal holography and topology change, when freed from the graph-theoretic combinatorics of 'spin networks',etc., revealing itself as being the 
Covering Theory of String/M Theory in AP Theory.Then some of the basic difficult problems of LQG (see,e.g. \cite{DT}) will be by-passed or solved intrinsically {\it ab initio}.

iii) Similarly AP-holography,our gauge-gravity correspondence,is so sharp and general,it can also be considered to be a particle-wave,particle-field 'duality': the discrete Artin presentation $r$,an 'extended' particle,i.e. a crystallic,non-topological 'quantum string',determines the smooth 4D spacetime universe $W^4(r)$. We can say:

{\it The 'particle'/'quantum string' $r$ is a hologram of the 4D spacetime universe $W^4(r)$; compare again to \cite{M2},p.63.}

This makes the fundamental AP construction of $W^4(r)$ from $r$,a rigorous background independent analogue of 'The Wave Function of the Universe',
\cite{HH},without using path integrals.
These analogies and the ones in section 3 already support many of the fundamental concepts of String/M theory as well as QCD,(see \cite{Wi3},p.1577,\cite{Wil1}, \cite{M4}), in a conceptually very simple,rigorous purely mathematical way:these concepts are {\it "here to stay"},\cite{L},p.9; the concepts of radical holography,universal topological change,etc. have to be present in any Quantum Geometry/ Gravity theory.

iv) We consider $h(r):\Omega_n\to \Omega_n$  to be the AP-analogue (where the rigorous Cayley-von Dyck intersection procedure is substituted for the non-rigorous Feynman summation procedure)  of 't Hooft's 'planar dominant Feynman diagrams' on the sphere $S^2$ (see \cite{Wi3},p.1577,\cite{M2},p.8, Ashtekar's remark on p.10 of \cite{Au}).
Thus AP Theory in its abstract,but rigorous,universal way,realizes 't Hooft's 'bold' Conjecture of
relating 4D Quantum Gauge Theory to String Theory,\cite{Wi3},p.1577, \cite{Wi10},p.25,\cite{Sm},p.44,\cite{S3}.

v) It is interesting to observe (compare to,e.g.,\cite{Mo}) that the truth of the 3D Poincar\'e  Conjecture shows that
there is no analogue of 4D  Black Hole singularities in AP Theory; although it is obvious from our Hilbert Vth-like construction,that 4D smooth singularities are avoided,they could have
'perversely' re-appeared as follows: if the Poincar\'e Conjecture were false,i.e.,if there existed an Artin Presentation $r$ presenting the trivial group,but such that the boundary
of the corresponding 4D smooth manifold $W^4(r)$ were not homeomorphic to the 3-sphere $S^3$,then this fact would imply that this smooth $W^4(r)$ would have a serious,unremovable
singularity: one would not be able to 'close' this smooth 4-manifold {\it smoothly}.See also \cite{Wi15}.

vi) AP Theory has rigorous analogues to all the features indicated in Witten's Fig.1 d) on p.25 of \cite{Wi11}. 

All the above make AP Theory a strong candidate for contributing to:{\it  "..the core geometrical ideas which underlie string theory, the way Riemannian geometry underlies general relativity"},see \cite{Wi12}.

In fact,perhaps it is not an exageration to call AP Theory,due to its discrete,group-theoretic conceptual simplicity,an "Erlanger Program" for Quantum Geometry and Quantum Gravity.It emerges as a new geometry {\it "..one that aligns with the new physics of string theory; "}, see \cite{Gr},pp.231,232.

In AP Theory, pure discrete group theory,i.e. 'symmetry' in its purest form, has been 'maximized' metamathematically; compare to \cite{Wil5},\cite{Wil4}.

{\it Metamathematically it is the most basic, simplest,'outermost' enveloping {\it framing} of Quantum Gravity}.

This basic intrinsic, autonomous mathematical consistency should also, a fortiori, encompass the logic, if rigorous, of any {\it empirical} physical evidence.\cite{FRS},p.24.

AP Theory in a way solves Atiyah's "joker in the pack" mystery of why low-dimensional manifold theory should be relevant to modern physics,\cite{A3},p.15: the, for each $n$, {\it infinitely generated} graded groups of topology changing transitions and interactions in AP Theory, take the place of the classical infinite dimensional Hilbert space of Quantum Theory.

The 4D metamathematical 'sporadicity' of AP Theory is much more universal than the important known ones: $Spin(4)$=$SU(2)\times SU(2)$,\cite{DK},p.7, or that the  braid group $B_n$ has a non-trivial amalgamation only when $n=3,4$.\cite{KPS}.

In our main section 3 we point out more rigorous specific analogies with the Clay YM problem and how these might impede its solution as stated and/or lessen its importance in modern physics.

We thank J.S.Calcut and C.E.Hough for very helpful conversations.

\section{The pure Mathematics of AP Theory}

In this section,we list,reference and explain some of the main,purely mathematical, but physically relevant, facts of
AP Theory,that have been rigorously proven in the refereed literature, (\cite{W1},\cite{CW},\cite{C2},..)

In a purely algebraic sense,AP Theory starts off as a sub-theory of the very basic,discrete
theory of Finitely Presented Groups,which begins with the concept of a presentation (see,e.g.,\cite{KMS},chapter 1) of a discrete group $G$ ,whose meaning is the following:

Let $F_n$ denote the free group on  the $n$ generators $x_1,\dots,x_n$ and let $w_1,\dots,w_m$  be $m$ words in $F_n$ ;let $N$ be the normal subgroup of $F_n$ ,which is the intersection of all normal subgroups of $F_n$  which contain all the $w_1\dots, w_m$;
then one says $<x_1,\dots,x_n\vert w_1\dots,w_m>$ presents $G$, is called a {\it presentation}  of $G$, if the factor group $F_n/N=G$. This is the Cayley-von Dyck process.[KMS,p.12]

It is important to notice here that a presentation of the trivial group,i.e. the case when it happens that $N =F_n$, can be,a priori, as complicated as a presentation of any arbitrary group  and that the concept of infinity is used here, when saying:"intersection of {\it all} normal subgroups".

Evidently in Group Theory,Presentation Theory is more basic, canonical and intrinsic than Representation Theory.

If $m=n$ a presentation $r=<x_1,\dots,x_n\vert r_1,\dots,r_n>$ is called an {\it Artin presentation},if,in $F_n$, the following group-theoretic equation holds:$$x_1\dots x_n=r_1^{-1}x_1r_1\dots r_n^{-1}x_nr_n$$.

As already realized by Artin himself, {\it it is an equation in the free group $F_n$ that actually  defines and characterizes pure,framed,(colored) braids.} See also \cite{MS}.

See \cite{W1},\cite{W2},\cite{W3},\cite{CW},\cite{C1},\cite{C2}, \cite{C3}, \cite{C4}, \cite{Ar}, and \cite{R},Appendix, for many examples.

The set of Artin presentations on $n$ generators is denoted by $R_n$, the group so presented by $\pi(r)$; $A(r)$ denotes the $n\times n$ integer matrix obtained by abelianizing the presentation $r$.

$A(r)$ is always symmetric and  determines the integer quadratic form of the compact 4-manifold $W^4(r)$; every symmetric, integer $n\times n$ matrix is an $A(r)$, where $r$ lies in $R_n$; $A(r)$ determines the $Z$-homology of both $W^4(r)$ and its boundary,$M^3(r)$.

We consider the integer, symmetric matrix $A(r)$ to be the analogue of the  Hilbert space binary forms (of QM) in AP Theory; it is all that remains of them under the radical reductivity of AP Theory,compare to \cite{Wi17},p.9,\cite{Wi2},p.4.

We say $r$ is 'a Torelli',if $A(r)$ is the zero matrix.

The $r$ in $R_n$ can be multiplied in a very non-trivial way,$r\cdot r'$ again being an Artin presentation, and so that $R_n$ is canonically isomorphic to $P_n\times Z^n$, where  $P_n$ denotes the pure braid group on $n$ strands (see \cite{W1},p.227]).The Torelli form a subgroup isomorphic to the (infinitely generated if $n>2$) commutator subgroup of $P_n$ and it is indeed a subgroup of the classical Torelli group of a 2D closed,orientable surface of genus $n$,hence the name. Multiplication of $r$ by a Torelli does not change $A(r)$ and hence preserves the Z-homology of both $W^4(r)$ and $M^3(r)$.

The fact that AP Theory is discrete and is graded by the positive integers,i.e. is 'cone-like',allows one to use mathematical induction,to stabilize,and to use computer based methods.See \cite{W1},\cite{W2}, \cite{W3}, \cite{CW}, \cite{C1}, \cite{C2}, \cite{C3}, \cite{C4}, for many examples of such computations.

Recall that a group is called perfect,if its abelianization is trivial.
From work of Milnor,see\cite{W1},p.227, we have the following 'Triality' fact:
{\it If the group $\pi(r)$ is finite and perfect,then it is either trivial,or isomorphic to I(120),the binary icosahedral group,i.e. the fundamental group of Poincar\'e's Z-homology 3-sphere}. 

The smooth 4D topological part of AP Theory starts as explained in the Introduction
(see \cite{W1},\cite{CW},\cite{C1}, \cite{C2},for the rigorous technical mathematical details) and we add the following in order to stress its topological 3D and 4D importance: 

i) Let $M^3(r)$ denote the boundary of  $W^4(r)$,the smooth 4-manifold determined by the Artin
presentation $r$ via the 2D diffeomorphism $h(r):\Omega_n\to \Omega_n$;it is always a connected,orientable,closed 3-manifold and every such 3-manifold can be so obtained.\cite{W1}, \cite{W3}, \cite{R},Appendix.To posess an Artin presentation characterizes the fundamental groups of such $M^3$ and the Artin presentation {\it actually determines the 3-manifold} up to diffeomorphism,not just its fundamental group. 

In particular,the theory of closed,orientable 3-manifolds,is strictly speaking,{\it not} an autonomous 3D theory,since at the same time of its definition,the smooth 4D manifold $W^4(r)$ is defined, amalgamating  this 3D theory with the much more physically relevant theory of 4D smooth manifolds.Thus the Hamilton-Thurston-Perelman program  aquires more physical importance,and hence perhaps a simpler solution.

ii)The $n+1$ boundary components of $\Omega_n$ define a link $L(r)$ of $n+1$ knots, $k_i(r)$,($i$=0,1,2,..,$n$),in $M^3(r)$; Gonz\'alez-Acu\~na showed (see \cite{C1}),{\it given any link $L$ in any closed,orientable 3D manifold $M^3$,there exists an Artin presentation $r$ such that $M^3$ = $M^3(r)$ and $L$ is a sublink of $L(r)$. In particular,any knot in any $M^3$ can be obtained this way}.

iii) When  $A(r)$ is unimodular,then $M^3(r)$ is a Z-homology 3-sphere,and the groups and peripheral structures of the above knots,have very simple,computer friendly presentations, \cite{W1}, p.226,  (where framings do not have to be 'put in by hand'),avoiding self-linking problems and heeding the admonitions of Penrose,Witten, et al.against symmetry-destroying 'skein methods'in Knot and Linking theory,when it is used in physics.

iv) The Torelli, denote one by $t$, by transforming $W^4(r)$,$M^3(r)$ into $W^4(t\cdot r)$,$M^3(t\cdot r)$, provide a very general,unpredictable and subtle theory of {\it topology changing} transitions, and form the analogue of Morse Theory in AP Theory, which is much sharper than that of any known $3+1$ TQFT.

Multiplying by a Torelli always preserves the $Z$-homology of $W^4(r)$ and $M^3(r)$ but usually changes the topology of $W^4(r)$ and $M^3(r)$,and the knots $k_i(r)$ in $M^3(r)$,but they can also just change certain things and leave others invariant:

{\bf Example:}Consider $s\in R_4$ and the Torelli $t\in R_4$ given by $s_1=(x_1x_3)^2x_2s_2$, $s_2=(x_1(x_2x_3)^2x_2^2)^{-1}$, $s_3=(x_2x_3x_2)^{-1}x_4s_4$,  $s_4=x_4^{-2}x_2x_3x_2(x_2x_3)^{-2}$, and $t_1=(x_4^{-1},x_1(x_1x_2x_3)^{-1})$,$t_2=(x_1,x_4)=t_3$, $t_4=(x_1^{-1},x_4x_1x_2x_3)t_3$; (here $(x,y)=x^{-1}y^{-1}xy$) and we use the computer algebra system MAGMA to do the group-theoretic computations).

MAGMA shows that $\pi(s) =1$, $M^3(s)=S^3$ and that  all the knot groups of the $k_i(s)$ are isomorphic to $Z$,except that of  $k_3(s)$,  which is isomorphic to that of the trefoil in $S^3$.

MAGMA,and some simplification by hand, gives $r=t\cdot s$ as:$$r_1=(x_4^{-1},x_1(x_1x_2x_3)^{-1})(x_4,x_1)(x_2x_3)^2x_2r_2,$$ $$r_2=(x_2x_3x_2^2)^{-1}((x_1^{-1},x_4x_1x_2x_3),x_4^{-1})(x_4^{-1},x_1)(x_1x_2x_3)^{-1},$$ $$r_3=(x_2x_3x_2)^{-1}(x_4x_1x_2x_3,x_1^{-1})x_4r_4,$$ $$r_4=x_4^{-2}(x_1^{-1},x_4x_1x_2x_3)x_2x_3x_2(x_2x_3)^{-2}(x_1,x_4);$$

Now we again have $\pi(r)=1$ and $M^3(r)=S^3$, however the (non-amphicheiral) trefoil $k_3(s)$ of $M^3(s)$ has been transformed by the Torelli $t$ to a (amphicheiral) figure-8 knot $k_3(r)$ in $M^3(r)$; all the other knots stay trivial.

v) Another very important,physically relevant property of the Torelli transitions is the following:

{\it One can change the smooth structure of a smooth 4-manifold, but leave the underlying topological structure intact. The discrete pure group theory of AP Theory has the energy and power to juggle different 4D smooth structures on the same underlying 4D topological manifold}.\cite{C1},\cite{C2}.

This phenomenon should be considered to be the last vestige (in the radical reductivity of AP Theory) of any hypothetical gravitational Schroedinger Wave Equation,where now 4D smooth structures, 'powered' by the Torelli transitions,  are the analogues of gravitational waves. This also seems to solve the so-called Hierarchy Problem ("Why gravity is so weak") in AP Theory.

vi) The symmetric integer matrix $A(r)$ determines the quadratic form of $W^4(r)$,in particular its integer homology,as well as the integer homology of $M^3(r)$.Any integer quadratic form can be so obtained,which implies the very non-trivial fact that AP Theory has a discrete purely group-theoretic analogue of the fundamental,physically relevant Donaldson Theorem,\cite{DK},\cite{GS},\cite{FM},\cite{Wi4}, \cite{Wi7}, \cite{Wi8}, \cite{Wi13}, {\it despite the fact that there are no moduli spaces nor deRham theory} in AP Theory:

{\bf THEOREM \cite{W1},p.240, \cite{R},p.621}:{\it If A(r) is a symmetric,integer,unimodular matrix,prevented by Donaldson's theorem from representing the quadratic form of closed,smooth,simply-connected 4-manifold,then the group $\pi(r)$ is non-trivial; in fact,it has a non-trivial representation into the Lie Group SU(2).} 

There exist (necessarilly non-Artin)  presentations $w$ of the trivial group, where $A(w)$ = $E_8$,(see p.11 of \cite{C4}) and hence it is the Artin Equation above, in the free group $F_n$,  for the presentation $r$, that makes this theorem true.  

Thus, in particular, the discreteness of AP Theory is related to the complex numbers in a very non-trivial way.

The 'modularity' hinted at by this {\it  purely group-theoretic} theorem should be related to that of Borcherds,\cite{B} and that of Tomita-Takesaki theory as in Algebraic QFT, \cite{S2},\cite{Sum}.

We remark,that at this point in time, (although,e.g. the Casson invariant, can be described in purely AP-theoretic fashion,\cite{Ar}, \cite{C3},\cite{C4}), we have no purely AP-theoretic proof of this 'Langlands-ian' theorem,which relates the  Number Theory of Integer Quadratic Forms with Group Representations of the $\pi(r)$  into $SU(2)$; compare to \cite{GW}, \cite{W4}.  We still need the actual classic  analytic field-theoretic methods for this.\cite{T},\cite{DS}.

On the other hand, this also shows that the non-analytic,non-local smooth 4D Hilbert Vth Problem-like construction above,
is so subtle , sharp and metamathematically 'local',that certain important analytic,{\it field-theoretic}, physically relevant results,pertaining to a single, given smooth 4-manifold, are still present
and are actually {\it detected} by the discrete AP Theory.

This  goes well beyond the mere ad hoc Tietze-like methods of the Kirby Calculus.

Furthermore,(see \cite{CW}),any complex elliptic surface,e,g., the Kummer surface, can be obtained smoothly as a $W^4(r)$,with  boundary $S^3$,thus also proving the existence in AP Theory
of a discrete, purely group-theoretic analogue of Donaldson/ Seiberg-Witten {\it invariants}.\cite{FM},p.4, \cite{DK},p.376,
\cite{Wi4},\cite{Wi6},p.375.

Finally,although in this paper we will not resort to it,we note that a very non-trivial Covering Theory exists in AP Theory: the covering and lifting theory of the 2D diffeomorphisms $h(r):\Omega_n\to \Omega_n$,{\it which is non-trivial even if the group $\pi(r)$ is trivial}.In particular,the covering theory of the so-called 'class surface',i.e. the regular covering corresponding to the commutator subgroup of  the fundamental group of $\Omega_n$, and that corresponding to the normal closure in $F_n$ of the $r_i$, should be specially relevant here for obtaining rigorous proofs in String theory and LQG.

This makes AP Theory also a very sophisticated mathematical theory indeed.

\section{Intrinsic 4D Quantum YM Existence and Mass Gap}

In sections 1,2,we exhibited the existence of a mathematically rigorous, sporadic,purely group-theoretic,smooth $3+1$ theory,which unlike the classical gauge-theoretic approach (which uses the space of connections on a particular,fixed manifold) is,a priori, exterior,(as cobordism theory is) ,autonomous,intrinsic,graded by the positive integers,and  which is not related to any particular Lie group in the usual gauge-theoretic sense.

As mentioned in the introduction, due to the radical,universal holography described above,this theory can not,a priori,be described by a conventional analytic,infinitesimal Field Theory.

This is the AP-analogue to {\it 'producing a mathematically complete example of quantum gauge field theory in four dimensional space time'}, \cite{JW},p.5, i.e. the first part of the Clay YM Existence and Mass Gap problem.

Despite that our '4D Quantum Gauge Theory' is not {\it quantitative},we can still ask the other fundamental question:

{\it Do there exist natural qualitative analogues to the quantitative mass gap condition of the Clay YM Millenium problem?}

We quote \cite{JW},p.3:{\it "..it must have a "mass gap":namely there must be some constant $\Delta > 0$ such that every excitation of the vacuum has energy at least $\Delta$".}

Furthermore,if such qualitative analogues of this quantitative mass gap do exist,how do they affect the rigorous solution,if it exists,of the actual {\it quantitative} Clay mass gap problem as stated?

We start with the questions of Ashtekar,\cite{As},p.174,regarding Quantum Geometry:{\it "What are the atoms of geometry? What are the fundamental excitations?"}.

We show that by considering the $h(r):\Omega_n\to \Omega_n$ as the analogue of 'vacuum fluctuations/excitations' and their generation of the 4D smooth manifolds $W^4(r)$ as a 'giving gravitational mass' Higgs-like phenomenon,we obtain some of the most desired and important hypothetical {\it qualitative} consequences of the {\it quantitative} Clay mass gap, if it were true.

In other words,the fundamental Hilbert Vth Problem-like construction, from the topology-less Artin presentation $r$ (i.e.$r$ is of the {\it discrete purely group-theoretic} 'vacuum' with zero 'mass'), gives non-zero 'gravitational mass',i.e.4D smooth structures to the $W^4(r)$ in a universal Higgs-like way.(Compare to the heuristic quantitative arguments in \cite{K} and references therein).

{\it In AP Theory,'mass gap' is just its radical,vacuum based,universal holography dressed in physical jargon}.

{\bf In AP Theory,holography and mass gap are defined in unison;they are the same mathematical phenomenon.}

We point out some AP-analogies with the most important {\it qualitative} consequences of a positive solution to the actual quantitative Clay YM problem:

I.Our mass gap immediately gives a sharper version of the 'clustering property' of (\cite{JW},p.6,or \cite{Wi10},p.125) {\it "of the principle of exponential decay of correlations at long distances that makes it possible to deduce global results about four manifolds from a knowledge how the theory behaves on $R^4$. ..the mass gap is closely related to the behaviour of the Donaldson invariants on algebraic surfaces"}.
(see also \cite{Wi13},p.291,\cite{Wi6}).
In AP Theory,global,non-local 4D  results are deduced not from $R^4$,{\it but already,more holographically,from $R^2$},the plane and Pure Braid Theory,via the planar $h(r):\Omega_n\to \Omega_n$.
As pointed out above,the results of \cite{CW} show that there even  exists a non-trivial,purely
AP-theoretic theory of Donaldson/Seiberg-Witten invariants.

II.The above AP-analogy (in iv) of the Introduction) with 't Hooft's 'bold' conjecture which according to \cite{Wi10},p.25,{\it ".. if valid,it might give an effective way to demonstrate the mass gap"} and which {\it "..seems like much the most plausible known approach to the problem,but an answer along these lines is not yet in sight,even at a heuristic level"}.

III.However the strongest analogy is with the  most important desideratum of the Clay problem: namely relating the mass gap with Quantum Chromodynamics (QCD),i.e. 4D $SU(N)$ Quantum Gauge Theory,\cite{Wi3}, p.1577, ('YM Theory without SUSY') and its important properties such as {\it confinement} and {\it asymptotic freedom}. \cite{Wil1},\cite{Wil2},\cite{H1},\cite{H2}.

First we note that the $h(r):\Omega_n\to \Omega_n$,our 'vacuum fluctuations/excitations' are not only sophisticated purely mathematically,but also {\it physically}:

Since $h(r)$ is determined only up to isotopies of $\Omega_n$,keeping it fixed as the identity on the boundary,we obtain a {\it topological} analogue of {\it uncertainty} in AP Theory (instead of minimal length 'uncertainty')
for our vacuum fluctuations.Compare to \cite{AY},\cite{SW}.

Thus,unlike as usual in quantum physics,'uncertainty' in AP Theory is {\it deduced} from the vacuum fluctuations,not {\it used} to 'prove' their existence.

In AP Theory, holography and 'uncertainty' are related; compare to \cite{SW}.

If we dare call the generic fixed points of $h(r):\Omega_n\to \Omega_n$, 'quarks',(compare to Wilczek,\cite{Wil3},Kondo,\cite{K} and Susskind's 'partons' in \cite{S}), we immediately obtain a topological analogue of 'confinement':

{\it Although generically these fixed points do not dissapear under isotopies,they can not be individually determined due to the uncertainty above}.

This is a topologically more sophisticated explanation of {\it confinement} ('quarks can not be individually determined'),than the classical string-theoretic one:that an open string has to have two inseparable end points,i.e.quarks.\cite{Wi3},p.1577.

This should also explain,why the phenomenon of confinement resists being proven analytically by field-theoretic methods.\cite{Wil4},p.7,\cite{S1},pp.40-41.

Similarly,if we also dare call our 4D smooth manifolds $W^4(r)$,'gluons',then we have a topological resemblance to {\it asymptotic freedom},('that in very high energy reactions quarks and gluons interact very weakly'):

If we iterate $h(r)$ it is natural to suppose that the mathematical relations between the fixed point theory (Nielsen,Thurston,..) of $h(r)^m$ and the 4D smooth topology of the corresponding $W^4(r^m)$ will grow weaker.This hints at the existence of an abstract intrinsic 'non-linear Fourier transform',\cite{A3},p.14,and should be compared to Taubes, \cite{T},p.367,Wilczek, \cite{Wil3}, Kondo,\cite{K}, and \cite{GW}, \cite{W4}.

Relating the the generic fixed points of $h(r):\Omega_n\to \Omega_n$ (i.e.'quarks') to the smooth structures of the 4D manifolds $W^4(r)$ is a sporadic {\it 4D smooth topological},i.e. gravitational, version of the original Yang-Mills program of relating particles to the differential geometry of connections.

It seems natural to conjecture that AP Theory and QCD ,due to their conceptual simplicity,are related ,in the sense that AP Theory has rigorous mathematical  features that the still hypothetical 'highest temperature QCD' should have.\cite{Wil1},p.25; see also \cite{We},pp.13,14.

In conclusion,we have exhibited the existence of an absolute,intrinsic '4D Quantum YM Theory' with a qualitative analogue of 'mass gap',with all the physical analogies above. Due to the universality and rigorous mathematical conceptual simplicity of this theory,perhaps the actual quantitative Clay YM problem should be substituted by the more general one of the existence of a $3+1$ axiomatic,'constructive' QFT. \cite{J}, \cite{Ri}, \cite{FRS}, \cite{S2},\cite{Wi14},\cite{We}.

Due to the crucial mathematical fact that AP Theory is not a model (in particular, does not introduce any new axioms) and its radical holography and strong Torelli transitions and interactions, which bypass difficult UV problems with pure group theory and which,a priori,does not seem to mix well with,e.g., the Wightman axioms,(see \cite{S2},p.10),the fate of a mathematically rigorous axiomatic $3+1$ QFT is a serious and  more pressing problem in physics,than the more particular Clay YM problem.It seems reasonable to conjecture that any such $3+1$ QFT would have to be 'modular' with respect to the groups $\pi(r)$ and therefore its rigorous construction, if indeed it exists, would be difficult.

For the Geometric Langlands Theory corresponding to Intrinsic 4D Quantum YM Theory, instead of $N=4$ Super YM Theory,  see \cite{W4}.

\end{document}